\documentclass[11pt]{amsart}
\usepackage{verbatim, latexsym, amssymb, amsmath}
\usepackage{epsfig}
\def\R{\mathbb R}
\def\N{\mathbb N}
\def\Z{\mathbb Z}
\def\C{\mathbb C}

\def\H{\mathcal H}
\def\L{\mathcal L}
\def\x{\mathbf x}

\def\k{\mathbf k}
\def\vol{\mathrm{vol}}
\def\re{\mathrm{Re}}

\def\a{\,d\mathcal{H}^n}
\def\c{\,d\mathcal{H}^2}

\newtheorem*{thma}{Theorem A}
\newtheorem*{thmb}{Theorem B}
\newtheorem*{thmc}{Theorem C}

\newtheorem{thm}{Theorem}[section]
\newtheorem{lemm}[thm]{Lemma}

\newtheorem{prop}[thm]{Proposition}
\theoremstyle{remark}

\theoremstyle{definition}

\title{Singularities of Lagrangian Mean Curvature Flow: Monotone case}

\begin{document}

\author{{Andr\'e Neves} ${}^{\dagger}$}
\email{aneves@math.princeton.edu} 
\address{Fine Hall, Princeton University, Princeton, NJ 08544, USA}

\thanks{\quad\ ${}^{\dagger}$\ The author was partially supported by NSF grant DMS-06-04164.}

\pagestyle{headings}

\begin{abstract}

We study the formation of singularities for the mean curvature flow of monotone Lagrangians in $\C^n$. More
precisely, we show that if singularities happen before a critical time then the tangent flow can be decomposed
into a finite union of area-minimizing Lagrangian cones (Slag cones). When $n=2$, we can improve this result by
showing that connected components of the rescaled flow converge to an area-minimizing cone, as opposed  to
possible non-area minimizing union of Slag cones. In the last section, we give specific examples for which such
singularity formation occurs.

\end{abstract}
\maketitle \markboth{Singularities for Lagrangian mean curvature flow: the monotone case} {Andr\'e Neves}

\section{Introduction}

We study the formation of singularities for Lagrangian mean curvature flow when the initial condition is a
Lagrangian in $\C^n$ for which the Liouville class equals the Maslov-class. The standard examples of such
Lagrangians are the Clifford torus in $\C^2$ which, under mean curvature flow, shrink to a point at time 1/2. In
this paper, we show that if the singularity occurs for some time $T<1/2$, then the tangent flow at the
singularity has the same properties as the tangent flow at a singularity when the initial condition is a
zero-Maslov class Lagrangian (see \cite{neves}). More precisely, we show that the tangent flow decomposes into a
finite union of Slag cones (Theorem A) and, when $n=2$, we improve this result and show that {\em connected}
components of the rescaled flow converge to a single SLag cone (Theorem B) as apposed to a union of SLag cones
which might have different Lagrangian angles and hence not be area-minimizing. If this property did not hold, it
would be hopeless to expect any kind of regularity theory for the Lagrangian at the singular time $T$. A similar
result was obtained in \cite{neves} in the context of rational Lagrangians having oscillation of the Lagrangian
angle smaller than $\pi$.

In the last section, we give examples of Lagrangians which are Hamiltonian isotopic to a Clifford torus and for
which the type of singularities described above occurs (Theorem C). These examples were communicated and
presented to Groh, Schwarz, and Smoczyk, in Leipzig 2005 and they have recently been included in their joint
work \cite{GSSZ}. In that paper, the authors show that Type I singularities do not occur for the flow of
monotone Lagrangians before a certain predetermined time. They also assume some rotational symmetry for the
initial condition and do a detailed study of the flow in this setting.

The organization of the paper is as follows. In the next section we introduce notation and derive some evolution
equations. In Section \ref{statement} we state Theorem A, Theorem B, and Theorem C. Section \ref{A} is devoted
to the proof of Theorem A and in Section \ref{B} we prove Theorem B. Theorem C is proven in  Section
\ref{examples}.

The author would like to thank Mu-Tao Wang for many interesting discussions.

\section{Preliminaires}
The standard complex structure and standard symplectic form  on $\C^n$ are denoted by $J$ and $\omega$
respectively. We consider also the closed complex-valued $n$-form given by
$$\Omega\equiv dz_1\wedge\ldots\wedge dz_n$$ and the Liouville form given by
$$
\lambda\equiv\sum_{i=1}^{n}x_idy_i-y_idx_i, \quad d\lambda=2\omega,
$$where $z_j=x_j+iy_j$ are complex coordinates for $\C^n$.

 A smooth $n$-dimensional submanifold $L$ in $\C^n$ is called {\em Lagrangian}
if $\omega_L=0$ and this implies that (see, for instance, \cite{HL})
$$\Omega_L=\exp(i\theta)\vol_L,$$
where $\vol_L$ denotes the volume form of $L$ and the multivalued function $\theta$ is called the {\em
Lagrangian angle}. Nevertheless, $d\theta$ is a well defined closed form on $L$ and its cohomology class is
denoted by {\em Maslov class}. If $\theta\equiv\theta_0$, then $L$ is calibrated by
$$\re\,\left( \exp({-i\theta_0})\Omega\right)$$ and hence area-minimizing. In this case, $L$ is referred as being
{\em Special Lagrangian}.

Note that the  Lagrangian condition implies that $\lambda$ is a closed form on $L$. A Lagrangian $L$ is said to
be {\em monotone} if
$$[\lambda]=c[d\theta]$$ for some positive constant $c$. When $c=1$ we call $L$ a {\em normalized monotone Lagrangian} and
note that every monotone Lagrangian can be rescaled in order to become normalized. The standard examples of
 monotone Lagrangians in $\C^2$ are the Clifford torus
 $$T_r\equiv \{(z_1,z_2)\in\C\,\mid\,|z_1|=r\mbox{ and }|z_2|=r\}.$$
The Clifford torus is normalized monotone when $r=1$.

Let $L_0$ be a  monotone Lagrangian and denote by $(L_t)_{0\leq t<T}$ a solution to Lagrangian mean curvature
flow. Denoting by $F_t$ the normal deformation by mean curvature, we have
\begin{multline*}
\frac{d}{dt}\int_{F_t(\gamma)}\lambda  = \frac{d}{dt}\int_{\gamma}F_t^*\lambda=\int_{\gamma}\L_H F_t^*\lambda
  =\int_{\gamma}dF_t^*(H\lrcorner\lambda)+F_t^*(H\lrcorner 2\omega)
  \\ = -2\int_{\gamma}dF_t^*\theta_t=-2[d\theta_t]
\end{multline*}
for every $[\gamma]$ in $H_1(L_t)$. Therefore, the monotone condition is preserved under the flow and we can
combine this with the fact that $[d\theta_t]$ is in $H^1(L_t,2\pi\Z)$ in order to obtain
\begin{lemm}\label{smoczyk}
For all $t<T$,
\begin{equation*}
[\lambda]=(c-2t)[d\theta_t] \quad\mbox{in }H^1(L_t).
\end{equation*}
\end{lemm}
 This
computation was communicated to the author by Smoczyk and can be seen as a higher codimension analogue of the
fact that curve shortening flow decreases enclosed area linearly. We obtain from Lemma \ref{smoczyk} the
existence of smooth functions $\gamma_t$ on $L_t$ such that
\begin{equation}\label{f_t}
d\gamma_t=\lambda+(2t-c)d\theta_t
\end{equation}
or, equivalently,
$$\nabla \gamma_t= (Jx)^{\top}+\nabla \theta_t$$
where $(Jx)^{\top}$ denotes the tangential projection on $T_xL_t$.
 The following computation can also be found in \cite{GSSZ}.

\begin{prop}\label{evol}
The functions $\gamma_t$ can be chosen so that
$$
\frac{d\gamma_t}{dt}=\Delta \gamma_t.
$$
\end{prop}

\begin{proof}
Assume without loss of generality that the family of functions $(\gamma_t)$ is smooth with respect to the time
parameter. Denoting by $d^*$ the negative adjoint of $d$, it is know that $d \theta_t$ satisfies the equation
(see \cite{smo1})
$$\frac{d (d\theta_t)}{dt}= d d^* (d\theta_t).$$
Furthermore, we can see from \cite[Lemma 6.2]{neves} that
$$d^*\lambda=H\lrcorner\lambda,$$
from which it follows that
\begin{align*}
\frac{d \lambda}{dt} & = \L_H \lambda = d(H\lrcorner\lambda)+H\lrcorner 2\omega
  =  d d^*\lambda-2d\theta_t.
\end{align*}
Combining this with identity \eqref{f_t}, we obtain that
\begin{align*}
    \frac{d( d\gamma_t)}{dt}= d d^*(\lambda+(2s-c)d\theta_t)=d\Delta \gamma_t
\end{align*}
and thus, we can add a time dependent constant to each $\gamma_t$ so that the desired result follows.

\end{proof}

\section{Statement of results}\label{statement}
If we consider a solution to mean curvature flow $(L_t)_{0\leq t <T}$ for which the initial condition is a
normalized monotone Lagrangian,  then Lemma \ref{smoczyk} shows that some sort of singular behavior happens at
time $t=1/2$. We focus in this paper on singularities that happen at the point $(x_0, T)$ in space-time with
$T<1/2$. The results we present next show that these type of singularities have the same behavior as finite time
singularities for the mean curvature flow of zero-Maslov class Lagrangians. We now give an heuristic argument
for why this phenomena is somehow expected.

Consider any point $x_1$ in $\C^n$ and define
$$\bar L_s\equiv \exp(s)(L_{(1-\exp(-2s))/2}-x_1).$$
Then, $\bar L_s$ solves
$$\frac{d x}{ds}=H+ \x^{\bot},$$
where $\x^{\bot}$ denotes the orthogonal projection on $(T_x \bar L_s)^{\bot}$ of the vector determined by the
point $x$ in $\C^n$. The change of basepoint does not change the cohomology class of $\lambda$ and so, we obtain
by Lemma \ref{smoczyk} that $\bar L_s$ is a normalized monotone Lagrangian. This implies the existence of smooth
functions $\bar\gamma_s$ such that, using Proposition \ref{evol},
$$J\nabla\bar\gamma_s=H+ \x^{\bot}\quad\mbox{and}\quad\frac{d\bar\gamma_s}{ds}=\Delta\bar\gamma_s+2\gamma_s.$$
This computations highlights the analogy between the family $(\bar L_s)$ and a solution to mean curvature flow
of zero-Maslov class Lagrangians. Note that $(\bar L_s)$ exist for all time precisely when $T=1/2$ and so,
arguing heuristically, finite-time singularities for zero-Maslov class Lagrangians ``correspond'' to
singularities at time $T<1/2$ for normalized monotone Lagrangians.

In order to study finite time singularities for mean curvature flow, take any sequence $(\sigma_j)$ going to
infinity and consider
\begin{equation}\label{Ljs}
L^{j}_s\equiv\sigma_j (L_{T+s/\sigma_j^2}-x_0) \quad \mbox{ for }-\sigma_j^2 T<s<0,
\end{equation}
which is still a solution to Lagrangian mean curvature flow (called the {\em rescaled flow}). For the rest of
this paper, the initial condition for the flow $L_0$ will have the property that, for some constant $C_0$,
$$
\H^n\bigl(L_0\cap B_R(0)\bigr)\leq C_0 R^n
$$
for all $R>0$

 Arguing informally, the following theorem states that a sequence of
rescaled flows at a singularity converges weakly to a finite union of integral Special Lagrangian cones.

\begin{thma} Let $L_0$ be a normalized  monotone Lagrangian in $\C^n$ for which the Lagrangian mean curvature flow develops a
singularity at time $T<1/2$. For any sequence of rescaled flows $(L^j_s)_{s<0}$ at a singularity, there exists a
finite set of Lagrangian angle $$\left\{\exp(i\bar\theta_1),\ldots,\exp(i\bar\theta_N)\right\}$$ and integral
Special Lagrangian cones
$$L_1,\ldots,L_N$$
such that, after passing to a subsequence, we have for every smooth function $\phi$ with compact support, every
$f$ in $C^2(S^1)$, and every $s<0$,
$$
\lim_{j \to \infty}\int_{L^j_s}f\left(\exp(i\theta_{j,s})\right)\phi\a=\sum_{k=1}^N m_k
f\left(\exp(i\bar\theta_k)\right)\mu_k(\phi)
$$
where $\mu_k$ and $m_k$ denote the Radon measure of the support of $L_k$ and its multiplicity respectively.

Furthermore, the set $\left\{\exp(i\bar\theta_1),\ldots,\exp(i\bar\theta_N)\right\}$ does not depend on the
sequence of rescalings chosen.
\end{thma}

When $n=2$ we can strengthen this result  and show that on each connected components of $L^j_s\cap B_R(0)$ the
Lagrangian angle converges to a constant. Arguing heuristically, this property shows that the formation of
singularities exhibits  a ``nice" behavior. A similar property was proven in \cite{neves} in the context of
rational Lagrangians in $\C^n$ with oscillation of the Lagrangian angle less than $\pi$.
\begin{thmb}
Let $L_0$ be a normalized  monotone Lagrangian in $\C^2$ for which the Lagrangian mean curvature flow develops a
singularity at time $T<1/2$. The following property holds for all $R>0$ and almost all $s <0$. \par For any
sequence $\Sigma^j$ of connected components of $B_{4R}(0)\cap L^j_s$ that intersect $B_R(0)$, there exists a
Special Lagrangian cone $\Sigma$ in $B_{2R}(0)$ with Lagrangian angle $\bar\theta$ such that, after passing to a
subsequence,
$$
\lim_{j \to \infty}\int_{\Sigma^j}f\left(\exp(i\theta_{j,s})\right)\phi\a=m
f\left(\exp(i\bar\theta)\right)\mu(\phi)
$$
for every $f$ in $C(S^1)$ and every smooth $\phi$ compactly supported in $B_{2R}(0)$, where $\mu$ and $m$ denote
the Radon measure of the support of $L$ and its multiplicity respectively.
\end{thmb}

Without the connectedness assumption on each $\Sigma_j$ the theorem would fail for trivial reasons. The
requirement that $\Sigma_j$ intersects $B_R(0)$ is necessary in order to prevent its area to converge to zero.

From the proof of Theorem $A$, it follows that for almost all $s<0$ and all $R>0$,
$$\lim_{j\to\infty}\int_{L^j_s\cap B_R(0)}|\nabla \theta_{j,s}|^2\c=0.$$
We remak that this alone is not enough to show Theorem B because the sequence of smooth manifolds $L^j_s$ are
becoming singular when $j$ goes to infinity and so, no Poincar\'e inequality holds with a constant independent
of $j$. This means we cannot conclude that, on each connected component of $L^i_s$, the Lagrangian angles
$\theta_{i,s}$ converge to a constant. To overcome this obstacle we will need to explore the fact that a
function $\gamma_t$ satisfying \eqref{f_t} is well defined.

Finally, assuming some rotational symmetry, we construct examples that show that the type of singularities
considered in Theorem A and Theorem B occur.

\begin{thmc}
There is a normalized monotone Lagrangian $L_0$ which is Hamiltonian isotopic to a Clifford torus and for which
the Lagrangian mean curvature flow develops a singularity at the origin at time $T<1/2$. The rescaled flow is
transverse union of two Lagrangian planes with the same Lagrangian angle.
\end{thmc}

\section{Proof of Compactness Theorem A}\label{A}

Before proving Theorem A we recall Huisken's monotonicity formula \cite{huisken} valid for any smooth family of
$k$-dimensional submanifolds $(N_t)_{t\geq 0}$ moving by mean curvature flow in $\R^m$. Consider the backward
heat kernel
$$\Phi_{x_0,T}(x,t)=\frac{1}{(4\pi(T-t))^{k/2}}e^{-\frac{|x-x_0|^2}{4(T-t)}}.$$
When $(x_0,T)=(0,0)$, we denote it simply by $\Phi$. The following formula holds
$$\frac{d}{dt}\int_{N_t} f_t \Phi_{x_0,T} = \int_{N_t}
\left(\frac{d}{dt} f_t -\Delta
  f_t-\left|H+\frac{(\x-\x_0)^{\bot}}{2(T-t)}\right|^2 f_t\right)\Phi_{x_0,T}$$
where $f_t$ is a smooth function with polynomial growth at infinity and $(\x-\x_0)^{\bot}$ denotes the
orthogonal projection on $(T_x N)^{\bot}$ of the vector determined by the point $(x-x_0)$ in $\R^m$.

The aim of this section is to prove Theorem A, which we restate it  for the sake of convenience.

\begin{thma} Let $L_0$ be a normalized  monotone Lagrangian in $\C^n$ for which the Lagrangian mean curvature flow develops a
singularity at time $T<1/2$. For any sequence of rescaled flows $(L^j_s)_{s<0}$ at a singularity, there exists a
finite set of Lagrangian angle $$\left\{\exp(i\bar\theta_1),\ldots,\exp(i\bar\theta_N)\right\}$$ and integral
Special Lagrangian cones
$$L_1,\ldots,L_N$$
such that, after passing to a subsequence, we have for every smooth function $\phi$ with compact support, every
$f$ in $C^2(S^1)$, and every $s<0$
$$
\lim_{j \to \infty}\int_{L^j_s}f\left(\exp(i\theta_{j,s})\right)\phi\a=\sum_{k=1}^N m_k
f\left(\exp(i\bar\theta_k)\right)\mu_k(\phi),
$$
where $\mu_k$ and $m_k$ denote the Radon measure of the support of $L_k$ and its multiplicity respectively.

Furthermore, the set $\left\{\exp(i\bar\theta_1),\ldots,\exp(i\bar\theta_N)\right\}$ does not depend on the
sequence of rescalings chosen.
\end{thma}

\begin{proof}
Set $u(x)\equiv\langle J x_0, x\rangle$, where $x_0$ is the point at which the singularity occurs. According to
Proposition \ref{evol}, we can  choose a family of functions $\gamma_t$ so that
$$d\gamma_t=\lambda+(2t-1)d\theta_t-du\quad\mbox{and}\quad \frac{d\gamma_t}{dt}=\Delta \gamma_t.$$
A simple computation shows that
\begin{multline*}
\frac{d}{dt}\int_{L_t} \gamma_{t}^2 \Phi_{x_0,T}\a \\= -\int_{L_{t}}
 \left(2\lvert\nabla\gamma_{t}\rvert^2+\lvert H-{(\x-\x_0)^{\bot}}/{2(T-t)}\rvert^2 \gamma_{t}^2\right)\Phi_{x_0,T}\a
\end{multline*}
and hence
\begin{equation*}
\lim_{t\to T}\int_{t}^{T}\int_{L_s}\lvert\nabla\gamma_{s}\rvert^2\Phi_{x_0,T}\a ds=0.
\end{equation*}
From Huisken's monotonicity formula we know that
\begin{equation*}\lim_{t\to
T}\int_{t}^{T}\int_{L_s}\left|H+(\x-\x_0)^{\bot}/{(2(T-t))}\right|^2\Phi_{x_0,T}\a ds=0
\end{equation*}
 and thus, because the identity
$$J\nabla \gamma_{t}=-(\x-\x_0)^{\bot}+(2t-1)H$$
implies that
$$|(2T-1)H|^2\leq 2|\nabla \gamma_t|^2+8(T-t)^2)|H+{(\x-\x_0)^{\bot}}/{(2(T-t))}|^2,$$
we obtain
\begin{equation}\label{luna}
 \lim_{t\to T}\int_{t}^{T}\int_{L_s}\left(|H|^2+|(\x-\x_0)^{\bot}/(T-t)|^2\right)\Phi_{x_0,T}\a ds=0.
\end{equation}

Fix any $a<b<0$. In terms of the sequence of flows $(L^j_s)$, it follows from the above identity that
\begin{equation}\label{ford}
\lim_{j\to\infty}2\int_{a}^{b}\int_{L^j_s}\left(|\x^{\bot}|^2+|H|^2\right)\Phi\a ds=0
\end{equation}
and hence, we can choose $a<0$ so that
$$
\lim_{i \to \infty}\int_{L^i_{a}\cap B_R(0)}\left(|\x^{\bot}|^2+|H|^2\right)\a=0
$$
for all positive $R$.

The rest of the proof follows from adapting the arguments used in the  proof of \cite[Theorem A]{neves}. The
first step consists in showing

\begin{prop}\label{general}
Let $(L^j)$ be a sequence of smooth Lagrangians in $\C^n$ such that, for some fixed $R>0$, the following
properties hold:
\begin{itemize}
\item[(a)] There exists a constant $D_0$ for which $$\H^n(L^j\cap B_{2R}(0))\leq D_0R^n$$ for all $j \in \N$;
\item[(b)]
$$\lim_{j \to \infty}\H^{n-1}\left(\partial L^j\cap B_{2R}(0)\right)=0$$
and $$\lim_{j \to \infty}\int_{L^j\cap B_{2R}(0)}|H|^2\a=0.$$
\end{itemize}
Then, there exist a finite set
$$\{\exp(i\bar\theta_1),\ldots,\exp(i\bar\theta_N)\}$$
and integral Special Lagrangians
$$L_1,\ldots,L_N$$
such that, after passing to a subsequence, we have for every smooth function $\phi$ compactly supported in
$B_R(0)$ and every $f$ in $C(S^1)$
$$
\lim_{j \to \infty}\int_{L^j}f(\exp(i\theta_{j}))\phi\a=\sum_{j=k}^N m_k f(\exp(i\bar\theta_k))\mu_k(\phi),
$$
where $\mu_k$ and $m_k$ denote, respectively, the Radon measure of the support of $L_k$ and its multiplicity.
\end{prop}
\begin{proof}
    A slightly different version of this proposition was proven in \cite[Proposition 5.1]{neves} where $L^j$ was
    assumed to be zero-Maslov class. The same proof applies straightforwardly provided we consider the sets
        $$\{x\in L^j\,\mid\,|\exp(i\theta_j(x))-\exp(i\bar\theta_1)|\leq\varepsilon\} $$
    and the functions $\cos\theta_j$ instead of the sets
        $$\{x\in L^j\,\mid\,|\theta_j(x)-\bar\theta_1|\leq\varepsilon\} $$
    and the functions $\theta_j$ respectively.
\end{proof}

A standard fact (see, for instance,  \cite[Lemma B.1]{neves}) implies the existence of a constant $D_0$ for
which
$$
\H^n\left(L^j_{a}\cap B_R(0)\right)\leq D_0 R^n
$$
for all positive $R$. Hence, we can apply Proposition \ref{general} to the sequence $(L^j_{a})$ and, after a
simple diagonalization argument, obtain a subsequence  (still indexed by $j$) for which there are integral
Special Lagrangians
$$L_1,\ldots,L_N$$
and a finite set
$$\{\exp(i\bar\theta_1),\ldots,\exp(i\bar\theta_N)\}$$
 such that for every smooth function $\phi$ compactly
supported  and every $f$ in $C^2(S^1)$
$$
\lim_{j \to \infty}\int_{L_a^j}f(\exp(i\theta_{j,a}))\phi\a=\sum_{j=k}^N m_k f(\exp(i\bar\theta_k))\mu_k(\phi),
$$
where $\mu_k$ and $m_k$ denote, respectively, the Radon measure and the multiplicity of $L_k$. The fact that
$$
\lim_{k \to \infty}\int_{L^j_{a}\cap B_R(0)}|\x^{\bot}|^2\a=0
$$
for all positive $R$ implies that the Special Lagrangians $L_k$ are all cones.

It is well known that (see, for instance, \cite{smo1})
$$
    \frac{d\exp(i\theta_{j,s})}{ds}=\Delta \exp(i\theta_{j,s})+\exp(i\theta_{j,s})\lvert H\rvert^2,
$$
and so, a simple computation shows that
$$P_{j,s}\equiv \frac{ df(\exp(i\theta_{j,s}))}{ds}-\Delta f(\exp(i\theta_{j,s}))$$
satisfies
\begin{equation}\label{pavement}
|P_{j,s}|+|\nabla f(\exp(i\theta_{j,s}))|^2\leq C|H|^2,
\end{equation}
 where $C=C(f)$. Hence,
    \begin{multline*}
        \frac{d}{ds}\int_{L^j_s}f(\exp(i\theta_{j,s}))\phi\a =- \int_{L^j_s}\langle\nabla f(\exp(i\theta_{j,s})),
        \nabla\phi\rangle\a\\+
        \int_{L^j_s}f(\exp(i\theta_{j,s}))\langle H,D\phi\rangle\a
        +\int_{L_s^j}P_{j,s}\phi\a\\-\int_{L^j_s}f(\exp(i\theta_{j,s}))|H|^2\phi\a.
    \end{multline*}
    It follows from \eqref{ford} and \eqref{pavement} that the integral in time of the right hand side converges to zero and
    so
\begin{multline*}
\lim_{j \to \infty}\int_{L^j_{b}}f(\exp(i\theta_{j,b}))\phi\a=\lim_{j \to
\infty}\int_{L^j_{a}}f(\exp(i\theta_{j,a}))\phi\a \\ =\sum_{k=1}^N m_k f(\exp(i\bar\theta_k))\mu_k(\phi)
\end{multline*}
for all  $b<0$.

We now prove the uniqueness statement of the theorem. Let
$$\bigl(\widehat L^j_s\bigr)_{s<0}$$ be another sequence of rescaled flows for which there are Special Lagrangian
cones
$$\widehat L_1,\ldots, \widehat L_P$$
and a finite set $\bigl\{\exp(i\hat\theta_1),\ldots,\exp(i\hat\theta_P)\bigr\}$ such that,  for every smooth
function $\phi$ with compact support, every $f$ in $C^2(S^1)$, and every $s<0$,
$$
\lim_{j \to \infty}\int_{\widehat L^j_s}f(\exp(i\theta_{k,s}))\phi\a=\sum_{k=1}^P \widehat m_k
f(\exp(i\hat\theta_k))\widehat\mu_k(\phi)
$$
where $\widehat\mu_k$ and $\widehat m_k$ denote the Radon measure of the support of $L_k$ and its multiplicity
respectively. We want to argue that
$$\bigl\{\exp(i\hat\theta_1),\ldots,\exp(i\hat\theta_P)\bigr\}=\{\exp(i\bar\theta_1),\ldots,\exp(i\bar\theta_N)\}.$$

For any real number $y$ and any integer $q$, we have the following evolution equation
$$\frac{d}{dt}(\cos\theta_t-y)^{q}=\Delta(\cos\theta_t-y)^{q}+U|H|^2,$$
where $$U\equiv q(\cos\theta_t-y)^{q-2}((\cos\theta_t-y)\cos\theta_t-(q-1)\sin^2\theta_t).$$
Therefore
\begin{multline*}
\frac{d}{dt}\int_{L_t}(\cos \theta_t-y)^q\Phi_{x_0,T}\a \\
=\int_{L_{t}} \left(U|H|^2-\lvert H-{(\x-\x_0)^{\bot}}/{2(T-t)}\rvert^2 \gamma_{t}^2\right)\Phi_{x_0,T}\a
\end{multline*}
  and we obtain
 from \eqref{luna} that
 $$\lim_{t\to T}\int_{L_t}(\cos \theta_t-y)^q\Phi_{x_0,T}\a$$
exists. Hence, from  scale invariance, we have for all $s,u<0$
\begin{multline*}
\lim_{j\to\infty}\int_{L^j_s}(\cos \theta^j_s-y)^q\Phi\a=\lim_{j\to\infty}\int_{\hat{L}^j_u}
(\cos{\theta}^j_u-y)^q\Phi\a\\=\lim_{t\to T}\int_{L_t}(\cos \theta_t-y)^q\Phi_{x_0,T}\a
\end{multline*}
and therefore, for all
positive integer $q$ and all $y$ in $\R$,
$$\sum_{k=1}^N m_k (\cos\bar\theta_k-y)^{q}\mu_k(\Phi)=\sum_{k=1}^P \widehat m_k
\bigl(\cos\hat\theta_k-y\bigr)^{q} \widehat\mu_k(\Phi).$$ Likewise, we can use the same arguments to show that
$$\sum_{k=1}^N m_k (\sin\bar\theta_k-y)^{q}\mu_k(\Phi)=\sum_{k=1}^P \widehat m_k
\bigl(\sin\hat\theta_k-y\bigr)^{q} \widehat\mu_k(\Phi)$$ for all positive integer $q$ and all $y$ in $\R$ and
this implies that
$$\bigl\{\exp(i\hat\theta_1),\ldots,\exp(i\hat\theta_P)\bigr\}=\{\exp(i\bar\theta_1),\ldots,\exp(i\bar\theta_N)\}.$$

\end{proof}

\section {Proof of Compactness Theorem B}\label{B}

This section is devoted to the proof of Theorem B.
\begin{thmb}
Let $L_0$ be a normalized  monotone Lagrangian in $\C^2$ for which the Lagrangian mean curvature flow develops a
singularity at time $T<1/2$. The following property holds for all $R>0$ and almost all $s <0$. \par For any
sequence $\Sigma^j$ of connected components of $B_{4R}(0)\cap L^j_s$ that intersect $B_R(0)$, there exists a
Special Lagrangian cone $\Sigma$ in $B_{2R}(0)$ with Lagrangian angle $\bar\theta$ such that, after passing to a
subsequence,
$$
\lim_{j \to \infty}\int_{\Sigma^j}f\left(\exp(i\theta_{j,s})\right)\phi\a=m
f\left(\exp(i\bar\theta)\right)\mu(\phi)
$$
for every $f$ in $C(S^1)$ and every smooth $\phi$ compactly supported in $B_{2R}(0)$, where $\mu$ and $m$ denote
the Radon measure of the support of $L$ and its multiplicity respectively.
\end{thmb}
\begin{proof}
    The fact that $L_0$ is monotone implies that, for each $L_s^j$,
        \begin{equation}\label{variacoes}
        [\lambda]=\sigma_j^2(1-2T+2s/\sigma_j^2)[d\theta_{j,s}].
        \end{equation}
    Moreover, identity \eqref{ford} implies that, after passing to a subsequence if necessary,
        \begin{equation}\label{coldplay}
        \lim_{j\to\infty}\int_{L^j_s\cap B_R(0)}|H|^2+|\x^{\bot}|^2\c=0
        \end{equation}
    for almost all $s<0$ and all $R>0$.

    The next proposition shows that connected components of  $B_{4R}(0)\cap L^j_s$ have a
    well defined Lagrangian angle and primitive for the Liouville form.

    \begin{prop}\label{smiths}
        For almost all $s<0$ and all $R>0$ the following property holds. Let $(\Sigma^j)$ be a sequence of connected components
        of $B_{4R}(0)\cap L^j_s$ that intersect $B_R(0)$. Then, for all $j$ sufficiently large, the Lagrangian angle
        $\theta_{j,s}$ is well defined, there is a smooth function $\beta_{j,s}$ so that $d\beta_{j,s}=\lambda$, and
            $$\H^{2}(\Sigma_j)\geq C_0 R^2,$$
        where $C_0$ is some universal constant.
    \end{prop}
    \begin{proof}
        We start by showing that an isoperimetric inequality  holds for almost all $s<0$.
        \begin{lemm}
            There is some universal constant $C$ such that, for all
        $j$ sufficiently large,
            $$\left(\H^2(A)\right)^{1/2}\leq C \H^{1}(\partial A),$$
        where $A$ is any open subset of $L^j_s\cap B_{6R}(0)$ with rectifiable boundary
        \end{lemm}
        \begin{proof}
        Recall that  according to the Michael-Simon Sobolev inequality there exists a
        universal constant $C$ so that
            $$\left(\H^2(A)\right)^{1/2}\leq C\int_A |H|\c +C\H^{1}(\partial A)$$
        and so
            $$\left(\H^2(A)\right)^{1/2}\leq C\left(\H^2(A)\right)^{1/2} \left(\int_A |H|^2\right)^{1/2}+C\H^{1}(\partial A). $$
        The claim follows from \eqref{coldplay}.
        \end{proof}
        Denote the intrinsic ball of radius $r$ around $x$ in $\Sigma_j$  by $\widehat{B}_{j}(x,r)$ and set
            $$\psi_j(r)\equiv \H^2\left(\widehat{B}_{j}(x,r)\right)$$
        which has, for almost all $r$, derivative given by
            $$\psi_j^{\prime}(r)=\H^{1}\left(\partial \widehat{B}_{j}(x,r)\right).$$
        Hence, we have for all $j$ sufficiently large, all $x$ in $\Sigma_j$, and all $r<R$
            $$(\psi_j(r))^{1/2}\leq C\psi_j^{\prime}(r),$$
        from which it follows that
            \begin{equation}\label{pixies}
                \H^2(\widehat{B}_{j}(x,r))\geq C_0 r^2,
            \end{equation}
        where $C_0$ is some universal constant. An immediate consequence is that $$\H^{2}(\Sigma_j)\geq C_0 R^2.$$

        The lower density bounds \eqref{pixies}  combined with the uniform upper bounds on $\H^2(\Sigma_j)$ imply the existence
        of some universal constant $k$ so that, for all $j$ sufficiently large,
            $$\mathrm{diam}\,\Sigma_j\leq kR.$$
        Hence, for all $j$ sufficiently large, we have that every loop $\gamma$ in $\Sigma_j$ satisfies
            $$\int_{\gamma}\lambda=0$$
        because
            $$\left|\int_{\gamma}\lambda\right|\leq 4R\,\mbox{length}(\gamma)$$
        and, due to \eqref{variacoes},
            $$\int_{\gamma}\lambda=\sigma_j^2(1-2T+2s/\sigma_j^2)2n\pi\quad\mbox{for some } n \in \Z.$$
        Note that from \eqref{variacoes}  we also obtain that
            $$\int_{\gamma}d\theta_{j,s}=0$$
        for every loop $\gamma$ in $\Sigma_j$ and thus the proposition is proven.
    \end{proof}

    Pick $s_0<0$ so that both \eqref{coldplay} and Proposition \ref{smiths} hold, and consider a sequence
     $\Sigma^j$ of connected components of $B_{4R}(0)\cap L^j_{s_0}$ that intersect $B_R(0)$.  We can apply Proposition
    \ref{general} and conclude that, after passing to a subsequence, $\Sigma_j$ converges weakly in $B_{2R}(0)$ to a
    finite union of Special Lagrangian cones that we denote by $\Sigma$. From Proposition \ref{smiths} we know that
    $\Sigma$ has positive measure. Next, using the argumentation in
\cite{neves}, we show that the support of $\Sigma$ is contained
    in the support of a Special Lagrangian cone.

    According to \eqref{variacoes}, we have
        $$[\lambda+2(s-s_0)d\theta_{j,s}]=(\sigma_j^2(2T-1)-2s_0)[d\theta_{j,s}]$$
    and thus, we can find a sequence $(b_j)$ converging to one so that $b_j^{-1}(\sigma_j^2(2T-1)-2s_0)$ is an
    integer, which implies that
        $$b_j^{-1}(\lambda+2(s-s_0)d\theta_{j,s})\in H^1(L^j_s,\Z) $$
    for all $s<0$. Hence, there is a multivalued function $\alpha_{j,s}$ so that
        $$d\alpha_{j,s}=b_j^{-1}(\lambda+2(s-s_0)d\theta_{j,s})$$ and
        $$u_{j,s}\equiv \cos \alpha_{j,s}$$
    is a well defined function.

     Assume without loss of generality that for $s=-1$ both \eqref{coldplay} and
    Proposition \ref{smiths} hold. As a result,  we have that, for all $R>0$ and all $j$ sufficiently large, both the
    Lagrangian angle and a primitive for the  Liouville are well defined on $L^j_{-1}\cap B_{2R}(0)$.
    \begin{lemm}\label{lemma1} There is a set
            $$
                \{(\cos\bar\beta_1,\sin\bar\beta_1) ,\ldots,(\cos\bar\beta_Q,\sin\bar\beta_Q)\}
            $$
        of distinct pairs and integral Special Lagrangian cones
            $$P_1,\ldots,P_Q$$ such that, after passing to a subsequence,
        we have for all smooth $\phi$ with compact support and all $f$ in $C(\R)$
    \begin{align*}
        & \lim_{j\to\infty}\int_{L^j_{-1}}f(\cos(b_j^{-1}\beta_{j,-1}))\phi\c=\sum_{k=1}^Q p_k f(\cos \bar\beta_k )\nu_{k}(\phi)\\
        & \lim_{j\to\infty}\int_{L^j_{-1}}f(\sin(b_j^{-1}\beta_{i,-1}))\phi\c=\sum_{k=1}^Q p_k f(\sin \bar\beta_k)\nu_{k}(\phi),
    \end{align*}
    where $\nu_k$ and the positive integer $p_k$ denote, respectively, the Radon measure of the support of $P_k$ and
    its multiplicity.
\end{lemm}
\begin{proof}
    Due to Proposition \ref{general}, Proposition \ref{smiths}, and identity \eqref{variacoes}, the proof given in \cite[Lemma 7.2]{neves}
    applies with no modifications.
\end{proof}

We can combine this lemma with Proposition \ref{general} and assume that, after a rearrangement of the supports
of the Special Lagrangian cones and its multiplicities, for all $\phi$ with compact support, all $f$ in $C(\R)$,
and all $y\in \R$,
 \begin{multline*}
    \lim_{j\to \infty}\int_{L^j_{-1}}f\left(\cos \left(b_j^{-1}({\beta_{j,-1}+2y\theta_{j,-1}})\right)\right)\phi\c\\=
    \sum_{k=1}^N m_k f(\cos(\bar\beta_k+2y\bar\theta_k))\mu_j(\phi),
 \end{multline*}
where $\mu_k$ denotes the Radon measure of the support of $L_k$ and the elements of the set
    $$\{(\cos\bar\beta_1,\sin\bar\beta_1,\bar\theta_1),\ldots,(\cos \bar\beta_N,\sin\bar\beta_Q,\bar\theta_N)\}$$
are all distinct.

    According to Proposition \ref{general}, the functions $\alpha_{j,s}$ can be chosen so that
        $$\frac{d u_{j,s}}{ds}=\Delta u_{j,s}+u_{j,s}|b_j^{-1}(\x^{\bot}+2(s_0-s)H)|^2.$$
    This evolution equation will be used to show

\begin{lemm}\label{lemma2}
    For all $\phi$ with compact support, all $f, h$ in $C^2(\R)$, and all $\bar\theta \in \R$,
    \begin{multline*}
        \lim_{j\to \infty}\int_{L^j_{s_0}}h(\cos(\theta_{j,s_0}-\bar\theta))f(\cos(b_j^{-1}\beta_{j,s_0}))\phi\c\\
        = \sum_{k=1}^N m_k
        h(\cos(\bar\theta_k-\bar\theta))f(\cos(\bar\beta_k-2(s_0+1)\bar\theta_k))\mu_k(\phi).
    \end{multline*}
\end{lemm}
\begin{proof}
    The proof is the same as in \cite[Lemma 7.3]{neves}. Set $$v_{j,s}\equiv h(\cos(\theta_{j,s}-\bar\theta))f(u_{j,s}).$$
    A simple computations shows that, for $\min\{s_0,-1\}\leq s<0,$
    $$P_{j,s}\equiv\frac{d v_{j,s}}{ds}-\Delta v_{j,s}$$
    is such that
    \begin{equation}\label{pjharvey}
        |P_{j,s}|+|\nabla v_{j,s}|^2\leq C(H^2+|x^{\bot}|^2),
    \end{equation}
    where $C=C(s_0,f,h)$.
    Hence,
    \begin{multline*}
        \frac{d}{ds}\int_{L^j_s}v_{j,s}\phi\c =- \int_{L^j_s}\langle\nabla v_{j,s},\nabla\phi\rangle\c+
        \int_{L^j_s}v_{j,s}\langle H,D\phi\rangle\c\\
        +\int_{L_s^j}P_{j,s}\phi\c-\int_{L^j_s}v_{j,s}|H|^2\phi\c.
    \end{multline*}
    Thus, it follows from \eqref{ford} and \eqref{pjharvey} that the integral in time of the right hand side converges to zero and
    so
    \begin{multline*}
        \lim_{j\to \infty}\int_{L^j_{s_0}}h(\cos(\theta_{j,s_0}-\bar\theta))f(\cos(b_j^{-1}\beta_{j,s_0}))\phi\c=
        \lim_{j\to\infty}\int_{L^j_{s_0}}v_{j,s_0}\phi\c\\
        =\lim_{j\to\infty}\int_{L^j_{-1}}v_{j,-1}\phi\c=\lim_{j\to\infty}\int_{L^j_{-1}}h(\cos(\theta_{-1,s}-\bar\theta))
        f(u_{j,-1})\phi\c\\
        =\sum_{k=1}^N m_k
        h(\cos(\bar\theta_k-\bar\theta))f(\cos(\bar\beta_k-2(s_0+1)\bar\theta_k))\mu_k(\phi).
    \end{multline*}
\end{proof}

Without loss of generality we can assume that
$$
    \cos(\bar\beta_1-2(s_0+1)\bar\theta_1),\ldots,\cos(\bar\beta_N-2(s_0+1)\bar\theta_N)
$$
are all distinct real numbers because this is true for all but countably many $s_0$. Moreover, we can use
\cite[Proposition A.1]{neves} and assume that on $\Sigma_j$ the sequence $(\cos (b_j^{-1}\beta_{j,s_0}))$
converges to a constant. Let $\gamma$ be the value of this constant and define $f \in C^2(\R)$ to be a
nonnegative cutoff function that is one in small neighborhood of $\gamma$ and zero everywhere else.

As a result, if $\mu_{\Sigma}$ denotes the Radon measure of $\Sigma$, we can use Lemma \ref{lemma2} to conclude
that, for every nonnegative test function $\phi$ with support in $B_{2R}(0)$,
\begin{multline*}
    \mu_{\Sigma}(\phi)=\lim_{j \to \infty}\int_{\Sigma_j}f(\cos(b_j^{-1}\beta_{j,s_0}))\phi\c
 \\\leq\lim_{j \to \infty}\int_{L^j_{s_0}}f(\cos(b_j^{-1}\beta_{j,s_0}))\phi\c\\
    =\sum_{k=1}^N m_kf(\cos(\bar\beta_k-2(s_0+1)\bar\theta_k))\mu_k(\phi).
\end{multline*}
Therefore, the fact that the support of $f$ can be made arbitrarily small and that the elements of the set
$$
\{\cos(\bar\beta_1-2(s_0+1)\bar\theta_1),\ldots,\cos(\bar\beta_N-2(s_0+1)\bar\theta_N)\}.
$$
are all distinct, implies that
$$\gamma=\cos(\bar\beta_{k_0}-2(s_0+1)\bar\theta_{k_0})$$ for a unique $k_0$. As a result,
$$
\mu_{\Sigma}(\phi)\leq m_{k_0} \mu_{k_0}(\phi)
$$
for all $\phi\geq 0$, which means that the support of $\Sigma$ is contained in the support of the Special
Lagrangian $L_{k_0}$.

In order to consider the proof completed we need to show convergence of the Lagrangian angle. Note that this
indeed requires an argument because the support of $L_{k_0}$ can be Special Lagrangian with  two different
orientations.

We know from \cite[Proposition A.1]{neves} that $(\cos (b_j^{-1}\beta_{j,s_0}))$ converges to a constant when
restricted to connected components of $L^j_{s_0}\cap B_{4R}(0)$ that intersect $B_{2R}(0)$.
 Let $\tilde\Sigma_j$ denote those connected components for which $(\cos (b_j^{-1}\beta_{j,s_0}))$ converges to
 $\gamma$. Using $\bar\theta=\bar\theta_{k_0}$
in Lemma \ref{lemma2} and keeping $f$ to be as defined above, we obtain that for every $\phi$ compactly
supported in $B_{2R}(0)$,
\begin{align*}
\lim_{j\to \infty}\int_{\tilde\Sigma_j}\cos(\theta_{j,s_0}-\bar\theta)\phi\c & =\lim_{j\to
\infty}\int_{L_{s_0}^j}\cos(\theta_{j,s_0}-\bar\theta_{k_0})f(u_{j,s_0})\phi\c\\
&=\lim_{j\to \infty}\int_{L_{-1}^j}\cos(\theta_{j,-1}-\bar\theta_{k_0})f(u_{j,-1})\phi\c\\
&=m_{k_0} \mu_{k_0}(\phi)=\lim_{j\to \infty}\int_{\tilde\Sigma_j}\phi\c.
\end{align*}
Hence, for all small $\varepsilon>0$,
$$\lim_{j\to\infty}\H^2(\{x\in\tilde\Sigma_{j}\,\mid\,\cos\theta_{j,s_0}(x-\bar\theta_{k_0})\leq 1-\varepsilon\})=0$$
and so the desired result follows.
\end{proof}

\section{Examples}\label{examples}

This section is devoted to the proof of

\begin{thmc}
There is a normalized monotone Lagrangian $L_0$ which is Hamiltonian isotopic to a Clifford torus and for which
the Lagrangian mean curvature flow develops a singularity at the origin at time $T<1/2$. The rescaled flow is
transverse union of two Lagrangian planes with the same Lagrangian angle.
\end{thmc}

\begin{proof}

Given an embedded curve $\gamma$ in the complex plane, it is easy to see that
$$
L=\{(\gamma\cos \alpha, \gamma \sin \alpha)\;|\; \alpha \in \R/2\pi\Z\}
$$
is a Lagrangian in $\C^2$ which will be embedded if $\gamma$ is preserved by the antipodal map. If we evolve $L$
by mean curvature flow, then each  $L_t$  is also rotational symmetric  and the corresponding $\gamma_t$ evolve
according to
\begin{equation}\label{flow}
\frac{d\gamma}{dt}=\k-\x^{\bot}/|x|^2,
\end{equation}
where $\k$ is the curvature of $\gamma$ and $\x^{\bot}$ denotes the orthogonal projection of the position vector
$x$ on the orthogonal complement of $T_x\gamma_t$. A detailed study of this equation is done in \cite{GSSZ}.
\begin{figure}\label{waits}
\centering {\epsfig{file=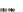, height=160pt}}\caption{$\bar L$ at the initial time and at the time of the
singularity.} \label{fig2}
\end{figure}

In \cite{neves} we showed that the curve in $\C$ given by
$$\bar\gamma\equiv\{s+i\,|\,s\in\R\}\cup \{s-i\,|\,s\in\R\}$$
(which gives rise to a zero-Maslov class Lagrangian $\bar L$) develops a finite time singularity at the origin
and its tangent flow is a union of two Lagrangian planes having the same Lagrangian angle (see Figure
\ref{fig2}).

Denote by $\gamma_{a}$ the family of ellipsoids
$$\gamma_a=\{x+iy\in\C\,|\, x^2/a^2+y^2/4\leq 1\}$$
and the correspondent family of Lagrangian torus by $L_a$. A simple calculation shows that, after choosing a
parametrization $\gamma_a (s)$,
$$\Omega_{L_a}=\frac{\gamma_a}{|\gamma_a|}\frac{\gamma_a^{\prime}}{|\gamma_a^{\prime}|}\vol_{L_a}\quad
\mbox{and}\quad\lambda_{L_a}=\langle i\gamma_a,\gamma_a^{\prime}\rangle ds.$$ Hence, one can easily check that
for all $a>0$, $L_a$ is monotone with the monotonicity constant $c_a$ equal to the area of the region enclosed
by $\gamma_a$ divided by $2\pi$. Denote the solution to Lagrangian mean curvature flow starting at $L_a$ by
$(L_{a,t})_{t\geq 0}$.

Choose $a$ large enough so that $c_a/2$ is strictly bigger than $\bar T$, the time at which the flow $(\bar
L_t)_{t\geq 0}$ with initial condition $\bar L$ develops a singularity. Note that  for all $t\leq \bar T$ and
while the solution exits, the region enclosed by $\gamma_{a,t}$ must contain the origin and, by Lemma
\ref{smoczyk}, its area is greater than a positive constant. The maximum principle for smooth solutions to
equation \eqref{flow} implies that $\bar\gamma_t$ cannot touch $\gamma_{a,t}$ while the solution remains smooth
and so $L_{a,t}$ must develop a singularity at time $T<c_a/2$ (see Figure \ref{fig3}).
\begin{figure}
\centering {\epsfig{file=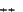, height=160pt}}\caption{$L_a$ at the initial time and at the time of the
singularity.} \label{fig3}
\end{figure}
The initial condition $L_a$ can be made normalized monotone after rescaling and in this case the singularity
will occur at time $T<1/2$. We are left to argue that the singularity happens at the origin. This follows from
either \cite[Theorem 1.14]{GSSZ} or the argument we present next.

For a short time, the curves $\gamma_{a,t}$ can be parameterized according to
$$
\gamma_{a,t}(s)=r_t(s)e^{is},\quad  s\in\R/2\pi\Z
$$
and, as can be seen from \cite[Lemma 4.6]{neves},
\begin{equation}\label{eqaux}
\frac{dr}{dt}=-\frac{\theta^{\prime}}{r}=\frac{r r^{\prime\prime}-2r^2-3(r^{\prime})^2}{r(r^{\prime})^2+r^3},
\end{equation}
where $\theta_t$ is the Lagrangian angle of $L_{a,t}$. Using the formula
$$\theta_0(s)=\arg(\gamma_t\gamma_t^{\prime})=2s+\arg(r_0^{\prime}+ir_0),$$
one can check that $\theta'_0(s)\geq 0$. Thus, arguing in the same way as in \cite[Lemma 5.2]{neves}, we obtain

\begin{lemm}\label{angle}While the solution to \eqref{eqaux} exists smoothly,
 $$\frac{dr}{dt}\leq0.$$
\end{lemm}

 If a singularity happens for equation \eqref{eqaux} at time $t_1$, then Lemma \ref{angle} implies that the
curve $\gamma_{t_1,a}$ can be parameterized as $r_{t_1}(s)e^{is}$. Thus, $t_1$ must also equal $T$ (the
singularity time \eqref{flow}) because otherwise $\gamma_{t_1,a}$ would be smooth and so $r_{t_1}(s)$ would also
be smooth.
\begin{lemm}\label{critical}
For any  $t<T$, $r_t(s)$ is nonincreasing when either $0\leq s\leq\pi/2$ or $\pi\leq s \leq 3\pi/2$, and
nondecreasing when either $\pi/2\leq s\leq\pi$ or $3\pi/2\leq s \leq 2\pi$.
\end{lemm}
\begin{proof}

Direct computation shows that, denoting $r_t^{\prime}$ by $u_t$,
$$
\frac{du_t}{dt}=\frac{u_t^{\prime\prime}}{(r^{\prime})^2+r^2}+u_t^{\prime}b(r_t,u_t,u_t^{\prime})+u_t
c(r_t,u_t,u_t^{\prime}),
$$ where the functions $b$ and $c$ are bounded for each $t<T$.

Because $\gamma_{t,a}$ is symmetric with respect to the coordinate axes, it is immediate to recognize that
$$u_t(0)=u_t(\pi/2)=u_t(\pi)=u_t(3\pi/2)=0$$
for all $t<T$. The result follows from the maximum principle.
\end{proof}

Suppose now that the singularity for the flow happens at a point $x_0=re^{i\alpha}$ with $r>0$. From Theorem A
we know that the tangent flow at the singularity is a union of planes and so, by White's regularity Theorem
\cite{white},
$$\limsup_{\delta\to 0}\frac{\H^1\left(\gamma_{T-\delta^2}\cap B_{\delta}(x_0)\right)}{2\delta}\geq 2.$$ This is impossible
because, using Lemmas \ref{angle} and \ref{critical}, we can argue like in \cite[Theorem 4.1]{neves} and
conclude that, for all $\delta$ sufficiently small and all $t<T$,
$$
\frac{\H^1\left(\gamma_{t}\cap B_{\delta}(x_0)\right)}{2\delta}\leq 3/2.
$$

 \end{proof}

\bibliographystyle{amsbook}

\vspace{20mm}

\end{document}